\documentclass[10pt]{amsart}
\usepackage{amssymb,amsfonts}
\usepackage{xypic}
\usepackage{amsthm}

\newtheorem{thm}{Theorem}[section]
\newtheorem{cor}[thm]{Corollary}
\newtheorem{pro}[thm]{Proposition}
\newtheorem{lem}[thm]{Lemma}

\theoremstyle{definition}
\newtheorem{defn}[thm]{Definition}

\newtheorem{exmp}[thm]{Example}

\newtheorem{rem}[thm]{Remark}

\newcommand{\g}{|G|}

\newcommand{\ind}{\text{ind}}
\newcommand{\res}{\text{res}}

\newcommand{\aaa}{\mathcal{A}}

\newcommand{\uu}{\mathcal{U}}

\newcommand{\zz}{\mathbb{Z}}

\newcommand{\thom}[1]{\text{Th} ( #1 )}
\newcommand{\euler}{\chi_{\text{Eu}}}
\newcommand{\rarr}{\rightarrow}
\newcommand{\col}{\colon}

\renewcommand{\contentsname}{Contents\\{\footnotesize\normalfont(A table
of contents should normally not be included)}}

\title{Generalized Artin and Brauer induction  \\ for compact Lie
groups}
\author{ Halvard Fausk }  \email{fausk@math.uio.no}
\address{Department of  Mathematics,  University of Oslo, 1053
Blindern, 0316 Oslo, Norway } 

\begin{document}

\begin{abstract} Let $G$ be a compact Lie group.
  We present two induction theorems for  certain generalized
  $G$-equivariant cohomology theories.
 The theory applies  to  $G$-equivariant $K$-theory $K_G$, and to
 the Borel cohomology associated
with any complex oriented cohomology theory. The coefficient ring of
$K_G$ is the representation ring $R(G)$ of $G$.
 When $G$ is a finite group the induction theorems for $K_G$
 coincide with
the classical Artin and Brauer induction theorems for $R(G)$.
\end{abstract}
\maketitle

\section{Introduction}
The Artin induction theorem, also called  Artin's theorem on induced characters, says that for any finite group $G$, the unit element in the representation ring $R(G)$, multiplied by the order of $G$,   is an integral  linear combination of elements induced  from $R(C)$, for cyclic subgroups $C$  of   $G$.  Similarly, the Brauer induction theorem, or Brauer's theorem on induced characters,  says that the unit element  in $R(G)$ is an integral linear
combination of elements induced  from $R(H)$, for subgroups  $H$ of $G$ which are extensions of  cyclic subgroups by $p$-groups.  These theorems in combination  with the double coset formula and the Frobenius reciprocity law  allow   $R(G)$ to be reconstructed integrally from $R(H)$, for subgroups $H$ which are extensions of cyclic groups by $p$-groups, and the restriction maps between them, and rationally from the $R(C)$, for
cyclic subgroups $C$, and the restriction maps between them. These reconstruction results are called restriction theorems
since they are  conveniently formulated using the restriction maps.

We present a generalization of the Artin and Brauer induction and restriction
theorems for the representation ring of a finite group $G$. The
generalization is in three directions. First, we give an induction
theory for a general class of equivariant cohomology theories; the
induction theorems apply to the cohomology groups of arbitrary
spectra, not just the coefficients of the cohomology theory. Second,
we extend the induction theory from finite groups to compact Lie
groups.
 Third, we allow induction from
more general classes of subgroups than the cyclic subgroups.
We use the following classes of abelian subgroups of $G$
characterized by the number of generators allowed: The class of the
maximal tori ($n =0$), and for each $n \geq 1$ the class of all
closed abelian subgroups $A$ of $G$ with finite index in its
normalizer and with a dense subgroup generated by $n$ or fewer
elements.

\

In section \ref{Liegroups} we collect some needed facts    on 
compact Lie groups. In section \ref{cohomology} we describe the
induction and restriction maps in homology and in cohomology.
 The induction theory  makes use of the  Burnside ring module
structure on equivariant cohomology theories. The Burnside ring is
isomorphic to the ring of homotopy classes of stable self maps of
the unit object $\Sigma^{\infty}_G S^0$ in the $G$-equivariant
stable homotopy category. In section \ref{Burnside} we recall some
alternative descriptions of the Burnside ring of a compact Lie group
$G$ and discuss some of their properties.

The following condition on a cohomology theory suffices to give the
induction theorems:
 We say that a ring spectrum $E$ has the
$0$-induction property if the unit map $\eta \col \Sigma^{\infty}_G
S^0 \rarr E$ pre-composed with a map $f\col \Sigma^{\infty}_G S^0
\rarr \Sigma^{\infty}_G S^0$ is null in the stable $G$-equivariant
homotopy category whenever the underlying nonequivariant map $f \col
\Sigma^{\infty} S^0 \rarr \Sigma^{\infty} S^0$ is null. We say that
a ring spectrum $E$ has the $n$-induction property, for some $n \geq
1$, if the unit map $\eta \col \Sigma^{\infty}_G S^0 \rarr E$
pre-composed with a map $f\col \Sigma^{\infty}_G S^0 \rarr
\Sigma^{\infty}_G S^0$ is null in the stable $G$-equivariant
homotopy category whenever the degree of $ f^A$ is $0$ for all
abelian subgroups $A$ of $G$, that is, the closure of a subgroup
generated by $n$ or fewer elements.
For example, singular Borel cohomology has the 0-induction property, equivariant $K$-theory has the  1-induction property, and the Borel cohomology associated with  a suitable height $n$-complex oriented ring spectrum (such as  $E_n$) has the  $n$-induction property.

 Let $E_G$ be a $G$-ring
spectrum satisfying the
 $n$-induction property.
   Pick  one  subgroup $A_i$
in each conjugacy class
 of the abelian
subgroups of $G$ with finite index in its normalizer and with a
dense subgroup generated by $n$ or fewer elements. Let $|G|_n$ be
the least common multiple of the order of the Weyl groups of these
abelian groups $\{ A_i \}$ (there are only finitely many such
subgroups by Corollary \ref{finite}). In section \ref{Artin} we
prove the following Artin induction theorem.
\begin{thm}  \label{introimage}
  The integer $\g_n$ times the unit element in
$E_{G}^{*}$ is in the image of the induction map \[ \textstyle\bigoplus\nolimits_i
\text{ind}_{A_i}^{G} \col \textstyle\bigoplus\nolimits_i E^{0}_{A_i} \rightarrow E^{0}_{G}
.\]
\end{thm}
Let $M_G$ be an $E_G$-module spectrum, and let $X$ be an arbitrary
 $G$-spectrum.
 There is a restriction map
\[ \res \col M^{\alpha}_G (X) \rightarrow \text{Eq} [
\textstyle\prod\nolimits_{i } M^{\res \, \alpha}_{A_i} (X)
\rightrightarrows \textstyle\prod\nolimits_{i,j ,g } M^{\res \,
\alpha}_{ A_i \cap g A_j g^{-1}} (X) ] . \] Here $\text{Eq} $
denotes the equalizer (the kernel of the difference of  the two parallel maps)  and $\alpha$ denotes the grading by a formal
difference of two finite dimensional real $G$-representations.
 The second product is over $i$, $j$ and over $g  \in G$.
  The maps in the equalizer are
the two restriction (composed with conjugation) maps. The Artin
induction theorem implies the following Artin restriction theorem.
\begin{thm} \label{Artinintroduction}  There exists a map
\[ \psi \col  \text{Eq} [ \textstyle\prod\nolimits_{i} M^{\res
\, \alpha}_{ A_i} (X)\rightrightarrows \textstyle\prod\nolimits_{i,j
,g } M^{\res \, \alpha}_{ A_i \cap g A_j g^{-1}} (X)] \rarr
M^{\alpha}_G (X)
\] such that both the composites $\res \circ \psi$ and $\psi
\circ \res$
are $\g_n$ times the identity map.
\end{thm}

The Brauer induction theorem is analogous to the Artin induction
theorem. At the expense of using a larger class, $\{ H_j \} $, of
subgroups of $G$ than those used for Artin induction, we get that
the unit element of $E^{\ast}_G $ is in the image of the induction
map from $\textstyle\bigoplus_j E^{\ast}_{H_j} $. As a consequence the
corresponding restriction map, $\res$, is an isomorphism. The exact
statements are given in section \ref{Brauer}.

The $G$-equivariant $K$-theory $K_G (X_+ )$ of a compact
$G$-CW-complex $X$ is the Grothendieck construction on the set of
isomorphism classes of finite dimensional complex $G$-bundles on
$X$. In particular, when $X$ is a point we get that $K_G (S^0 )$ is
isomorphic to the complex representation ring $R(G)$. The
induction and restriction maps for the equivariant cohomology theory
$K_G$ give the usual induction and restriction maps for the
representation ring $K_G (S^0 ) \cong R(G)$. The ring spectrum $K_G$
satisfies the 1-induction property. The resulting induction theorems
for $R(G)$ are the classical Artin and Brauer induction theorems.
The details of this example are given in section \ref{ktheory}.

 This work is
inspired by an Artin induction theorem used by Hopkins, Kuhn, and
Ravenel \cite{hkr}. They calculated the Borel cohomology associated
with certain complex oriented cohomology theories for finite abelian
groups; furthermore  they used an Artin restriction theorem to describe the
Borel cohomology, rationally, for general finite groups. We discuss
induction and restriction theorems for the Borel cohomology
associated with complex oriented cohomology theories when $G$ is a
compact Lie group in section \ref{Borel}. Singular Borel cohomology
is discussed in section \ref{singular}.

A Brauer induction theorem for the representation ring of a compact
Lie group was first given by G.~Segal \cite{seg1}. Induction
theories for $G$-equivariant cohomology theories, when $G$ is a
compact Lie group, have also been studied by G.~Lewis
\cite[sec.6]{lew}. He develops a Dress induction theory for Mackey
functors.

The idea to use the Burnside ring module structure to prove
induction theorems
  goes back to Conlon  and Solomon
\cite{con,sol} \cite[chapter 5]{ben}.

\section{Compact Lie groups} \label{Liegroups}
In this section we recall some facts about compact Lie groups and
provide a few new observations.
 We say that a  subgroup $H$ of $G$ is topologically
generated by $n$ elements (or fewer) if there is a dense subgroup of
$H$ generated by $n$ elements; e.g.~any torus is
topologically generated by one element. By a subgroup of a compact
Lie group $G$ we mean a closed subgroup of $G$ unless otherwise
stated. It is convenient to give the set of conjugacy classes of
(closed) subgroups of $G$ a topology \cite[5.6.1]{td}. Let $d$ be a
metric on $G$ so that the metric topology is equal to the topology
on $G$. We give the space of all closed subgroups of $G$ the
Hausdorff topology from the metric \[ d_H ( A , B ) = \sup_{a \in A
} \inf_{b \in B} d ( a , b) + \sup_{b \in B} \inf_{a \in A } d ( a ,
b) .\] The metric topology on the space of closed subgroups is
independent of the choice of metric.
 Let $\Psi G$ denote the space of conjugacy classes of closed
subgroups of $G$ given the quotient topology from the space of
closed subgroups of $G$. The space $\Psi G$ is a compact space with
a metric given by \[d_{\Psi} ( A , B ) = \inf_{g \in G} d_H ( A , g
B g^{-1} ) . \] The Weyl group $W_G H$ of a subgroup $H$ in $G$ is $
N_G H/H$. Let $\Phi G$ denote the subspace of $\Psi G$ consisting
 of conjugacy classes of
subgroups of $G$ with finite Weyl group. We have that $\Phi G$ is a
closed subspace of $\Psi G$ \cite[5.6.1]{td}.

 We denote the conjugacy
class of a subgroup $H$ in $G$ by $(H)$, leaving $G$ to be
understood from the context. Conjugacy classes of subgroups of $G$
form a partially ordered set; $(K) \leq (H)$ means that $K $ is
conjugate in $G$ to a subgroup of $H$. Note that a subgroup of a
compact Lie group $G$ cannot be conjugate to a proper subgroup of
itself. (There are no properly contained closed $n$-manifolds of a
closed connected $n$-manifold.)

A theorem of Montgomery and Zippin says that for any subgroup $H$ of
$G$ there is an open neighborhood $U$ of the identity element in $G$
such that all subgroups of $H U $ are subconjugate to $H$
\cite[II.5.6]{bre}\cite{mz}.

Let $K \leq H $ be subgroups of $G$. The normalizer $N_G K$ acts
from the left on $(G/H)^K$. Montgomery and Zippin's theorem implies
that the coset $(G/H)^K /N_G K $ is finite \cite[II.5.7]{bre}. In
particular, if $W_G K$ is finite, then $(G/H)^K$ is finite. The Weyl
group $W_G H$ acts freely on $ (G/H )^K$ from the right by $gH \cdot
nH = gnH$, where $gH \in (G/H)^K$ and $nH \in W_G H$. So $ |W_G H|$
divides $| (G/H )^K|$.
The following consequence of Montgomery and Zippin's theorem is
important for this paper. Let $G^{\circ}$ denote the unit component
of the group $G$.
\begin{lem}
\label{openpoint} Let $G$ be a compact Lie group.
 The  conjugacy class of any abelian
subgroup $A$ of $G$ with finite Weyl group is an open point in $\Phi
G$.
\end{lem}
\begin{proof}
 Fix a metric on $G$.
By Montgomery and Zippin's theorem there is an $\epsilon >0$ such
that if $K$ is a subgroup of $G$ and $d_{\Psi} ( (A) , (K) ) <
\epsilon $, then $K$ is
 conjugated in $G$ to a subgroup of $A$ that meets all the
 components of $A$.
Let $K$ be such a subgroup and assume in addition that it has finite
Weyl group. Then $K^{\circ} = A^{\circ}$ since $A < N_G K$ and $|W_G
K| $ is finite. Thus we have that $(K)=(A)$. Hence $(A)$ is an open
point in $\Phi G$.
\end{proof}

Since $\Phi G$ is compact we get the following.
\begin{cor} \label{finite}
There are only finitely many conjugacy classes of abelian subgroups
of $G$ with finite Weyl group.
\end{cor}
Given a subgroup $H$ of $G$ we can extend $H$ by tori until we get a
subgroup $K$ with $W_G K$ finite. This extension of $H$
 is unique up to conjugation. We denote the
conjugacy class by $\omega (H)$. The conjugacy class $\omega (H)$
does only depend on the conjugacy class of $H$. Hence we get a well
defined map $ \omega \col \Psi G \rarr \Phi G$. This map is
continuous \cite[1.2]{fol}. We say that $\omega (H)$ is the $G$
subgroup conjugacy class with finite Weyl group associated with $H$.
One can also show that the conjugacy class $\omega (H)$ is the
conjugacy class $(HT)$ where $T$ is a maximal torus in $ C_G H $
\cite[2.2]{fol}. This result implies the following.

\begin{lem} \label{omegaabelian}
The map $ \omega \col \Psi G \rarr \Phi G$ sends conjugacy classes
of abelian groups to conjugacy classes of abelian groups.
\end{lem}

We now define the classes of abelian groups used in the Artin
induction theory.

\begin{defn} \label{classes}
Let $\aaa G $ denote the set of all conjugacy classes of abelian
subgroups of $G$ with finite Weyl group. Let $\mathcal{A}_n G $
denote the set of conjugacy classes of abelian subgroups $A$ of $G$
that are topologically generated by $n$ or fewer elements and that
have a finite Weyl group. We let $\mathcal{A}_0 G $ be the conjugacy
class of the maximal torus in $G$.
\end{defn}  We often suppress $G$ from the
notation of $\aaa_n G$ and write $\aaa_n$. We have that $\aaa_n G =
\aaa G$ for some $n$ by Corollary \ref{finite}.
\begin{exmp} \label{cartan}  The  topologically cyclic subgroups
of $G$ are well understood. They were called Cartan subgroups and
studied by G.~Segal in \cite{seg1}. The following is a summary of
some of his results: All elements of $G$ are contained in a Cartan
subgroup. An element $g$ in $G$ is called regular if the closure of
the cyclic subgroup generated by $g$ has finite Weyl group. The
regular elements of $G$ are dense in $G$. Two regular elements in
the same component of $G$ generate conjugate Cartan subgroups. The
map $S \mapsto G^{\circ} S/ G^{\circ} $ gives a bijection between
conjugacy classes of Cartan subgroups and conjugacy classes of
cyclic subgroups of the group of components $G/G^{\circ}$. In
particular, if $G$ is connected, then the Cartan subgroups are
precisely the maximal tori. The order $|S/S^{\circ}|$ is divided by
$|S/G^{\circ}|$ and divides $|S/G^{\circ}|^2$ \cite[p.~117]{seg1}.
For example the nontrivial semidirect product $S^1 \rtimes \zz / 2
 $ has Cartan subgroups $S^1$ and (conjugates of) $ 0 \rtimes \zz
/ 2  $.
\end{exmp}

\begin{lem} \label{absplit}
If $A$ is a compact abelian Lie group, then it splits as
\[ A \cong A^{\circ} \times \pi_0 (A) .\]
\end{lem}
\begin{proof}
Since $A$ is compact we have that $\pi_0 (A) \cong \textstyle\bigoplus_i \zz /
p_{i}^{n_i} $. The unit component $A^{\circ}$ is a torus. We
construct an explicit splitting of $ A \rightarrow \pi_0 (A)$. Let
$a_i \in A$ be an element such that $a_i$ maps to a fixed generator
in $\zz / p_{i}^{n_i}$ and to zero in $\zz / p_{k}^{n_k}$ for all $k
\not= i$. Then $a_i$ raised to the $ {p_{i}^{n_i}}$ power maps to
zero in $\pi_0 (A)$,
hence is in the torus $ A^{\circ} $. There is an element $b_i \in
A^{\circ} $ such that
\[ a^{p_{i}^{n_i}}_i = b^{p_{i}^{n_i}}_i . \] Set $\bar{a}_i = a_i
b^{-1}_i$ and define the splitting $\pi_0 (A) \rightarrow A $ by
sending the fixed  generator of $ \zz / p_{i}^{n_i} $ in $ \pi_0 (A)$ to
$\bar{a}_i$. Since $A$ is commutative this gives a well defined
group homomorphism.
\end{proof}
The splitting in Lemma \ref{absplit} is not natural.
\begin{lem} \label{abelianfamily}
Let $G$ be a compact Lie group, and let $A \leq B$ be abelian
subgroups of $G$ such that $ \omega (B) $ is in $\aaa_n$. Then
$\omega ( A) $ is in $\aaa_n$.

\end{lem}  \begin{proof}
The minimal number of topological generators of an abelian group $A$
is equal to the minimal number of generators of the group of
components
 of $A$ by Lemma \ref{absplit}.
We assume without loss of generality that $B$ has finite Weyl group.
The unit component $B^{\circ}$ of $B$ is contained in the normalizer
of $A$. Hence the component group of a representative for the
conjugacy class $\omega (A)$ is isomorphic to a quotient of a
subgroup of $\pi_0 (B)$. The result then follows.
\end{proof}

We introduce several different orders for compact Lie groups. T.~tom
Dieck has proved that for any given compact Lie group $G$ there is
an integer $n_G$ so that the order of the group of components of the
Weyl group $W_G H$ is less than or equal to $n_G$ for all closed
subgroups $H$ of $G$
  \cite{tdf}.
\begin{defn} \label{order}
The order $|G|$ of a compact Lie group $G$ is the least common
multiple of the orders $ |W_G H |$ for all $(H) \in \Phi G$. For any
nonnegative integer $n$ let $\g_n$ be the least common multiple of $
|W_G A|$ for all $(A) \in \aaa_n G$.
\end{defn}
When $G$ is a finite group all these orders coincide and are equal
to the number of elements in $G$.
\begin{rem}
Let $T$ be a maximal torus in $G$. Then we have that
\[ N_G T /(G^{\circ} \cap  N_G T ) \cong G / G^{\circ} \]
since all maximal tori of $G$ are conjugated by elements in
$G^{\circ}$. Hence the number of components $| G / G^{\circ} | $ of
$G$ divides the smallest order $|G|_0 =|N_G T /T|$. The order
$|G|_m$ divides $|G|_n$ for $0 \leq m\leq n$.
\end{rem}
\begin{exmp} For compact Lie groups the various orders might be
different. An example is given by $SO (3)$. The only conjugacy
classes of abelian subgroups of $SO (3) $ with finite Weyl group are represented by the subgroups
$ \zz /2 \oplus \zz /2 $ and $S^1 $ of $SO (3)$. The normalizers of these
subgroups are $( \zz /2 \oplus \zz /2 ) \rtimes \Sigma_3$ and $ S^1
\rtimes \zz /2 $, respectively \cite[5.14]{td}. So the Weyl groups
have orders 6 and 2, respectively. Hence $|SO (3)|_n =2$ for $n = 0,
1 $ and $|SO (3)|_n =6$ for $n \geq 2 $. By taking cartesian
products of copies of $SO(3)$ we get a connected compact Lie group
with many different orders. The order $|SO (3)^{\times N} |_{2m}$ is
$2^N 3^m$ for $m \leq N$ and $6^{ N}$ for $m \geq N$. The abelian
subgroups of $SO (3)^{\times N}$ with finite Weyl group are product
subgroups obtained from all conjugates of $S^1$ and $ \zz /2 \oplus
\zz /2 $, and furthermore all subgroups of these product subgroups
so that each of the $n$ canonical projections to $SO (3) $ are
conjugate to either $S^1$ or $ \zz /2 \oplus \zz /2 $ in $SO (3) $.
\end{exmp}

\section{$G$-equivariant cohomology theories} \label{cohomology}
We work in the homotopy category of $G$-spectra indexed on a
complete $G$-universe. Most of the results  used are from
\cite{lms}.   If   $X$  is a based $G$-space the  suspension spectrum, $\Sigma^{\infty}_G X$, is simply denoted  by $X$.

We recall the definition of homology and cohomology theories
associated with a $G$-equivariant spectrum $M_G $. Let $X$ and $Y$ be
$G$-spectra. Let $ \{X , Y\}_G$ denote the group of stable (weak)
$G$-homotopy classes of maps from $X$ to $Y$.
We grade our theories by formal differences of $G$-representations.
For brevity let $\alpha $ denote the formal difference $V - W$ of
two finite dimensional real $G$-representations $V$ and $W$. Let $S^V$ denote the one point compactification of $V$.  Let $
S_{G}^{\alpha} $ denote the spectrum $ S^{-W} \Sigma_{G}^{\infty}
S^V $.
 The homology is \[ M^{G}_{\alpha} ( X) = \{ S_{G}^{\alpha} ,  M_G
\wedge X \}_G \cong \{ S_{G}^V , S_{G}^W \wedge M_G \wedge X \}_G .
\] The cohomology is
\[ M_{G}^{\alpha} ( X) = \{ S_{G}^{- \alpha} \wedge X ,
 M_G \}_G \cong  \{ S_{G}^W \wedge X , S_{G}^V \wedge M_G  \}_G .
\]

In this paper a ring spectrum $E$ is a spectrum together with a
multiplication $\mu \col E \wedge E \rarr E$ and a left unit $\eta
\col S_{G}\rarr E$ for the multiplication in the stable homotopy
category. We do not need to assume that $E$ is associative nor
commutative. An $E$-module spectrum $M$ is a spectrum with an action
$ E \wedge M \rarr M$ by $E$ that respects the unit and
multiplication. Let $E_G$ be a $G$-equivariant ring spectrum. The
coefficients $ E_{G}^{\alpha} = E^{G}_{ - \alpha} $ have a bilinear
multiplication that is $RO(G)$-graded and have a left unit element.
Let $M_{G}$ be an $E_G$-module spectrum. We have that $M_{G}^{\ast}
( X)$ is naturally an $E_{G}^{\ast} $-module, and $M^{G}_{\ast} (
X)$ is naturally an $E^{G}_{\ast} $-module for any $G$-spectrum $X$.

Let $M_G$ be a spectrum indexed on a $G$-universe $\uu$. For a
closed subgroup $H$ in $G$ let $M_H$ denote $M_G$ regarded as an
$H$-spectrum indexed on $\uu$ now considered as an $H$-universe. The
forgetful functor from $G$-spectra to $H$-spectra respects the smash
product. A
  complete $G$-universe $\uu$ is also a complete $H$-universe
  for all closed subgroups $H$ of $G$.
  For lack of a reference we include an argument proving
  this well-known result.

  Let $V$ be an $H$-representation.  The manifold $G \times V$ has a
   smooth and  free
  $H$-action  given by  $h \cdot (g , v) = (g h^{-1} , h v ) $. It
  also has   a smooth $G$-action by letting $G$ act from
  the right on $G$ in $ G \times V$.
The $H$-quotient $G \times_H V$ is a smooth $G$-manifold
\cite[VI.2.5]{bre}. Now consider the tangent $G$-representation $W$
at $ (1, 0) \in G \times_H V$. This can be arranged so that $W$ is
an orthogonal $G$-representation by using a $G$-invariant Riemannian
metric on $G$. Compactness of $H$ gives that $V$ is a summand of $W$
regarded as an $H$-representation.

    We have the following
isomorphisms for any $\alpha$ and any $G$-spectrum $X$
\cite[XVI.4]{ala}:
 \[  M^{H}_{\res_{H}^G \alpha} ( X) \cong \{ G/H_+ \wedge S_{G}^{
\alpha} , M_G \wedge X \}_G  ,
\]
\[ M_{H}^{\res_{H}^G \alpha} ( X) \cong \{ G/H_+ \wedge S_{G}^{- \alpha}
\wedge X , M_G \}_G . \] The forgetful functor from $G$-spectra to
$H$-spectra respects the smash product.

We now consider induction and restriction maps. The collapse map $c
\col G/H_+ \rarr S_{G}$ is the stable map associated with the $G$-map
that sends the disjoint basepoint $+$ to the basepoint $0$, and
$G/H$ to $1$ in $S^0 =\{ 0 ,1 \}$. Let $\tau \col S_{G}\rarr G/H_+ $
be the transfer map \cite[IV.2]{lms}. We recall a construction of
$\tau$ after Proposition \ref{frobenius}.

 There is an induction map natural in the $G$-spectra $X$ and $M_G$:
\[ \ind_{H}^G \col  M_{\res_{H}^G  \alpha }^H  (X) \rightarrow  M_{\alpha }^G
(X). \] It is defined by pre-composing with the transfer map $S_{G}
\stackrel{\tau }{\longrightarrow} G/H_+$ as follows:
\[  M_{ \res_{H}^G \alpha }^H  (X) \cong \{ G/H_+ \wedge S_{G}^{\alpha}
, M_G \wedge X \}_G \stackrel{\tau^{\ast} }{\longrightarrow}
M_{\alpha }^G (X) . \] There is a restriction map natural in the 
$G$-spectra $X$ and $M_G$:
\[ \res_{H}^{G} \col   M_{\alpha }^G  (X)  \rightarrow
 M_{\res_{H}^G  \alpha }^H  (X). \] It is defined by pre-composing
with the collapse map $G/H_+ \stackrel{c }{\rightarrow} S_{G}$. The
definition is analogous for cohomology. Alternatively, we can
describe the induction map in cohomology as follows:
\[ M^{\alpha }_G  (\tau \wedge 1_X ) \col
 M^{\alpha }_G  (G/H_+ \wedge X ) \rarr
M^{\alpha }_G ( X ) \] and the restriction map as
\[ M^{\alpha }_G  (c \wedge 1_X ) \col
 M^{\alpha }_G  (  X ) \rarr M^{\alpha }_G  (G/H_+ \wedge X ) \]
composed with the isomorphism \[ \{ S_{G}^{- \alpha} \wedge G/H_+
\wedge X , M_G \}_G \stackrel{(k\wedge 1_X)^* }{\longrightarrow} \{
G/H_+ \wedge S_{G}^{- \alpha} \wedge X , M_G \}_G \cong
M^{\res_{H}^{G} \alpha }_H ( X ) \] where $k \col G/H_+ \wedge
S^{-\alpha}_G \cong S^{-\alpha}_G \wedge G/H_+ $.

The classical Frobenius reciprocity law says that the induction map
$R(H) \rarr R(G)$ between representation rings is linear as an
$R(G)$-module, where $R(H)$ is given the $R(G)$-module structure via
the restriction map. In our more general context the Frobenius
reciprocity law says that the induction map $M^{\res_{H}^G \alpha}_H
( X) \rightarrow M^{\alpha}_G ( X)$ is linear as a map of $
E^{\ast}_G (S_{G}) $-modules (via the restriction map). We need the
following slightly different version.

\begin{pro} \label{frobenius} Let $M_G$ be a module over a ring
spectrum $E_G$. Let $ e \in E^{H}_{\res_{G}^{H} \alpha} $ and $ m
\in M^{G}_{\beta} ( X)$. Then we have that
\[ \ind_{H}^{G} (e) \cdot m =  \ind_{H}^{G} ( e \cdot \res_{H}^{G} m
) \] in $M^{G}_{\alpha + \beta} ( X)$.
The same result applies to cohomology.
\end{pro}
\begin{proof} Let $e \col   G/H_+ \wedge S^{ \alpha}_{G}
   \rightarrow E_{G}  $
represent the element $e$ in $E_{\res_{G}^{H}\alpha}^{H} $, and let
$ m \col
 S^{ \beta}_{G}  \rightarrow M_{G} \wedge X $ represent the
 element
 $m \in M^{G}_{\beta} (X)$.

We get that both products are
\[ S_{G}^{ \alpha} \wedge S_{G}^{ \beta} \stackrel{\tau \wedge 1
}{\longrightarrow}G/H_+ \wedge S_{G}^{ \alpha} \wedge S_{G}^{ \beta}
\stackrel{e \wedge m }{\longrightarrow } E_G \wedge M_G \wedge X
\] composed with the   $ E_G \wedge M_G \rarr M_G$.
The proof for cohomology is similar.
\end{proof}

We now describe the induction and restriction maps for homology
theories in more detail. This is used in section \ref{Borel}. We
have that
\[ \{ S_{G}^{\alpha} \wedge G/H_+ , E \wedge X \}_G \cong \{ S_{G}^{\alpha} ,
D(G/H_+ ) \wedge E \wedge X \}_G \] where $D (G/H_+) $ is the
Spanier-Whitehead dual of $G/H_+$. Using the equivalences $
S_{G}^{-\alpha} \wedge G/H_+ \cong G/H_+ \wedge S_{G}^{-\alpha} $ and $ D
(G/H_+) \wedge E \cong E \wedge D (G/H_+)$ we get an isomorphism
\[ E^{H}_{\res_{H}^{G} \alpha } (X) \cong
E^{G}_{\alpha} ( D(G/H_+) \wedge X ) .\]
 Under this isomorphism the  induction map is given by
 $E^{G}_{\alpha} ( D (\tau) \wedge 1_X )$, and the restriction map as
$E^{G}_{\alpha} ( D (c) \wedge 1_X )$. In the rest of this section
we recall a description of the transfer map \cite[IV.2.3]{lms} and
the Spanier-Whitehead dual of the collapse and transfer maps
\cite[IV.2.4]{lms}. Let $M$ be a smooth compact manifold without
boundary. In our case $M = G/ H$. There is an embedding of $M$ into
some finite dimensional real $G$-representation $V$
\cite[VI.4.2]{bre}. The normal bundle $\nu M$ of $M$ in $V$ can be
embedded into an open neighborhood of $M$ in $V$ by the equivariant
tubular neighborhood theorem. The Thom construction $\thom{
\xi}$ of a bundle $\xi$ on a compact manifold is equivalent to the
one point compactification of the bundle $\xi$. We get a map
\[ t' \col  S^V \rarr \thom{  \nu M }\] by mapping everything outside
of the tubular neighborhood of $M$ to the point at infinity.

The Thom construction of the inclusion map $ \nu M \rarr \nu M
\oplus TM \cong V \times M $ gives $ s'\col \thom{ \nu M } \rarr
S^V \wedge M_+ $. Let the pretransfer $t\col S_{G} \rarr S^{-V }
\wedge \thom{ \nu M } $ be $ S^{-V } \wedge t' $ pre-composed
with $S_{G} \cong S^{- V } \wedge S^V$ and let $s \col S^{-V } \wedge
\thom{ \nu M } \rarr M_+$ be the composite of $ S^{-V} \wedge s'
$ with $S^{-V} \wedge S^V \wedge M_+ \cong M_+ $.
  The transfer map $\tau$ is defined to be the
composite map \[ s \circ t \col S_{G} \rarr M_+ .\] We now let $M$ be
the $G$-manifold $G/H$.
Atiyah duality gives that the Spanier-Whitehead dual of $G/H_+$ is
equivalent to $S^{-V} \thom{ \nu G/H } $. When $G$ is finite
this is just $G/H_+$ itself. The proof of the equivariant Atiyah
duality theorem \cite[III.5.2]{lms} gives that $ D( c ) \simeq t$.
It is easy to see that $ D (s) \simeq s$. Hence we get
\[ D (\tau ) \simeq c \circ s \ \ \text{and}  \ \ D (c ) \simeq t .\]
 The discussion above gives the following.
 \begin{lem}  \label{homology}
 Let $c \col  G / H_+ \rarr S_{G}$,
 $ s\col  S^{-V} \wedge \thom{  \nu G/H } \rarr
 G/H_+ $, and  $t \col  S_{G} \rarr
S^{-V} \wedge \thom{ \nu G/H }$ be as above. Then
the restriction map in homology is
\[( t\wedge 1_X )_{\alpha}  \col  E^{G}_{\alpha} (X) \rarr
E^{G}_{\alpha} ( S^{-V} \thom{ \nu G/H } \wedge X )
\] and the induction map in homology is
\[ ( (c\circ s ) \wedge 1_X )_{\alpha} \col  E_{  \alpha}^G
( S^{-V} \thom{\nu G/H} \wedge X ) \rarr E_{\alpha}^G ( X)
\] composed with the isomorphism $ E_{  \alpha}^G
( S^{-V} \thom{\nu G/H } \wedge X ) \cong E^{H}_{res_{G}^{H}
\alpha} (X)$.
\end{lem}
If $G$ is a finite group, then $s$ is the identity map. Hence the
induction map is the induced map from the collapse map $c$, and the
restriction map is the induced map from the transfer map $ \tau$,
composed with the isomorphism $ E_{ \alpha}^G (G/H_+ \wedge X )
\cong E^{H}_{res_{G}^{H} \alpha} (X)$.

\section{The Burnside ring} \label{Burnside}
The stable homotopy classes of maps between two $G$-spectra are
naturally modules over the Burnside ring of $G$. We use this
Burnside ring module structure to prove our induction theorems.

We recall the following description of the Burnside ring $A (G) $ of
a compact Lie group $G$ from \cite{tdb,td,lms}. Let $ a (G) $ be the
semiring of isomorphism classes of compact $G$-CW-complexes with
disjoint union as sum, cartesian product as product, and the point
as the multiplicative unit object.
 Let $C( \Phi G ; \mathbb{Z} ) $, or
$C(G)$ for short, be the ring of continuous functions from the space
$\Phi G$ of conjugacy classes of closed subgroups of $G$ with finite
Weyl group to the integers $\zz$. Let $\euler (X)$ denote the Euler
characteristic of a space $X$. We define a semi-ring homomorphism
from $a(G)$ to $C(G)$ by sending $ X $ to the function $(H) \mapsto
\euler (X^H )$ \cite[5.6.4]{td}. This map extends to a ring
homomorphism ${\phi}' $ from the Grothendieck construction $b(G) $
of $a(G)$ to $C(G)$. The Burnside ring $A(G)$ is defined as $b(G) /
\ker {\phi}'$. We get an injective ring map $\phi \col A(G) \rarr
C(G)$. The image of $\phi$ is generated by $ \phi (G/H_+ ) $ for $H
\in \Phi G$. One can show that $C(G)$ is freely generated by $ |W_G
H |^{-1} \phi (G/H_+ ) $ for $H \in \Phi G$ \cite[V.2.11]{lms}.
 We have  that $ \g C (G) \subset
A(G)$ where $\g$ is the order of $G$ \cite[thm.2]{tdf}. We denote
the class $\chi (X)$ in $A(G)$ corresponding to a finite
$G$-CW-complex $X$ by $[X]$.
Define a ring homomorphism $d\col \pi^{G}_0 (S_{G} ) \rarr C(G)$
by sending
  a stable map $ f \col  S_{G} \rightarrow S_{G} $ to the function
that sends $(H)$ in $\Phi G$ to the degree of the fixed point map $
f^H$. It follows from Montgomery and Zippin's theorem that the maps
$\phi'$ and $d$ take values in continuous functions from $\Phi G$ to
$\zz$.

There is a map $\chi \col A(G) \rarr \pi^{G}_0 (S_{G} ) $ given by
sending the class of a compact $G$-CW-complex $X$ to the composite
of the transfer and the collapse map
\[ S_{G} \rightarrow  X_{+}
\rightarrow S_{G} . \] The map $\chi$ is the categorical Euler
characteristic \cite[V.1]{lms}. It turns out to be a ring
homomorphism \cite[V.1]{lms}. It has the property that the degree of
the $H$-fixed point of a map in the homotopy class $\chi (X)$ is
equal to $ \euler (X^H )$ (the ordinary Euler characteristic of the
fixed point space $X^H$) \cite[V.1.7]{lms}. A proof is given in
Lemma \ref{degree}. We have the following commutative triangle:
\[ \xymatrix{ A(G) \ar[rr]^{\chi} \ar[dr]_(.45){\phi} & &
\pi_{0}^G (S_{G} ) \ar[dl]^(.45)d \\ & C ( \Phi G ; \zz ) . & } 
\] A theorem,  due to Segal when $G$ is a finite group and to tom
Dieck when $G$ is a compact Lie group, says the map $\chi$ is an
isomorphism \cite[V.2.11]{lms}. This allow us to use the following
three different descriptions of elements in the Burnside ring:
\begin{enumerate}
\item Formal differences of
equivalence classes of compact $G$-CW-complexes. \item Stable
homotopy classes of self maps of $S_{G}$. \item Certain continuous
functions from $\Phi G$ to the integers.
\end{enumerate}

Since the Burnside ring $A(G)$ is isomorphic to $\{ S_{G} ,
S_{G} \}_G$ we have that $G$-equivariant cohomology and homology
theories naturally take values in the category of modules over
$A(G)$.

In the rest of this section we prove that the degree of the
$H$-fixed points of a stable map $ f \col S_{G} \rarr S_{G}$ is
the same as the degree of the $\omega (H)$-fixed point of $f$ for
any closed subgroup $H$ of $G$.
 The
following is well known.
\begin{lem} \label{abq} Let $X$ be any space with an action by a torus $T$.
If both $\euler (X) $ and $\euler (X^{T} )$ exist, then $\euler (X)
= \euler (X^{T} )$.
\end{lem} \begin{proof} Replace   $X$ by a weakly equivalent $T$-CW-complex.  We get that the quotient  complex $X / X^T$ is built
out of one single point $\ast$ and cells $ D^n \wedge T/A$ for a
proper subgroup $A$ of $T$. All nontrivial cosets of $T$ are tori
(of positive dimension). Hence all the cells have Euler
characteristic equal to 0 except for the point which has Euler
characteristic 1. The claim follows by the long exact sequence in
homology and the assumptions that both $\euler (X) $ and $ \euler
(X^{T} )$ exist.
\end{proof}

The following is a generalization of \cite[V.1.7]{lms}. They
consider closed subgroups with finite Weyl group.
\begin{lem} \label{degree} Let $X$ be a compact $G$-CW-complex and
let $f\col S_{G} \rarr S_{G} $ be a stable $G$-map in the stable
homotopy class $\chi (X) \in \pi^{G}_0 (S_{G} )$. Then
$\text{deg} ( f^L ) = \euler (X^L )$ for any closed subgroup $L$ of
$G$.
\end{lem}
\begin{proof}
   The  geometric fixed point functor $\Phi^L$
is a strong monoidal functor from the stable homotopy category of
$G$-spectra to the stable homotopy category of $W_G L$-spectra
\cite[II.9.12]{lms}. We also have that $ \Phi^L ( \Sigma^{\infty}_G
X_+ ) \cong \Sigma^{\infty}_{W_G L } X^{L}_+$ \cite[XVI.6]{ala}. The
forgetful functor from the stable $W_G H$-homotopy category to the
nonequivariant stable homotopy category is also strong monoidal. The
categorical Euler characteristic respects strong monoidal functors
\cite[3.2]{pic}. Hence we get that $ \deg ( f^L ) $ is equal to the
degree of the categorical Euler characteristic of the spectrum $
\Phi^L (\Sigma^{\infty}_G X )$ regarded as a nonequivariant
spectrum. This is $ \euler (X^L )$.
\end{proof}

\begin{pro} \label{omega}
Let $f\col S_{G} \rarr S_{G} $ be a stable $G$-map. Let $H$ be a
closed subgroup of $G$ and let $\omega (H)$ be the associated
 conjugacy class of subgroups with finite Weyl group. Then we have that
\[ \text{deg}  ( f^H )  = \text{deg} (f^{\omega (H) } ) .\]
\end{pro}
\begin{proof}  Let $K$ be a subgroup in the  conjugacy class
$\omega (H)$ so that $H<K$ and $K /H$ is a torus. Let $X$ and $Y$ be
finite $G$-CW-complexes such that $f$ is in the homotopy class $\chi
(X) - \chi (Y)$. Since
\[ \euler ( X^K ) = \euler ( (X^H)^{K/H} ) \text{ \ and \ }
\euler ( Y^K ) = \euler ( (Y^H)^{K/H} )
\] the  previous two lemmas give that \[
\text{deg} (f^K ) = \euler (X^K ) - \euler (Y^K) = \euler ( X^H ) -
\euler (Y^H) = \text{deg} ( f^H ) . \qedhere
\]
\end{proof}

We need the following corollary in section \ref{Borel}.

\begin{cor} \label{abeliandegree}
Let $f \col S_{G} \rarr S_{G} $ be a map such that $\deg ( f^A )
=0 $ for all abelian subgroups of $G$ with finite Weyl group. Then
  $f$ is null homotopic restricted to the $K$-equivariant
stable homotopy category for any abelian subgroup $K$ of $G$.
\end{cor} \label{null}
\begin{proof}   Proposition \ref{omega} together with Lemma
\ref{omegaabelian} give that the degree $ \deg ( f^A ) $ is $ 0 $
for all abelian subgroups $A$ of $G$. The claim follows since a self
map of $S^{0}_{K} $ is null homotopic if and only if  the degrees
of all its fixed point maps are 0 \cite[8.4.1]{td}.
\end{proof}

\section{Artin induction} \label{Artin}
We first introduce some conditions on ring spectra and then  prove
the Artin induction and restriction theorems. Recall Definition
\ref{classes}.
 Let $J_n$ be the $A(G)$-ideal  consisting of all
elements $\beta \in A(G)$ such that $\text{deg} (\beta^A ) =0 $ for
all $ (A) \in \aaa_n$. Let $J $ be the intersection of all $J_n$.

\begin{defn} \label{inductionproperty}
We say that a $G$-equivariant ring spectrum $E_{G}$ satisfies the
$n$-induction property if $ J_n E^{0}_G = 0$. We say that $E_G$
satisfies the induction property if $J E^{0}_G =0$.
\end{defn}
  Let  $\eta \col  S_{G} \rarr E_G$ be the unit
map of the ring spectrum $E_G$. Then  $E_G$ satisfies the
$n$-induction property if and only if the ideal $J_n$ is in the
kernel of the unit map
\[ \eta \col   A (G) \rarr E^{0}_G . \]
 If $E_G$ satisfies the $n$-induction property and
$E'_G$ is an $E_G$-algebra, then $E'_G$ also satisfies the
$n$-induction property.

 Let $e_{H} \col  \Phi G \rarr \mathbb{Z}$ be the function defined by
letting $e_H (H)=1$ and $e_H (K) =0 $ for $(K) \not= (H)$. We have
that $e_A$ is a continuous function for every $(A) \in \aaa$ since
$(A)$ is an open-closed point in $\Phi G$ by Lemma \ref{openpoint}.
Since $|G| C(G) \subset \phi A(G)$ we have that $|G| e_A \in \phi
A(G)$ for all $A \in \aaa$. When $G$ is a compact Lie group it turns
out that we can sharpen this result. Recall Definition \ref{order}.

\begin{pro} \label{idempotent}
Let $(K)$ be an element in $\aaa_n$. Then $\g_n e_{K}$ is an element
in $ \phi A(G)$. Moreover, the element can be written as
\[ \g_n  e_{K} = \textstyle\sum\limits_{(A) \leq (K)} c_A \phi (G/A) \]
where $ c_A \in \mathbb{Z}$ and $(A) \in \aaa_n$.
\end{pro}
\begin{proof}
Let $\mathcal{S} (K)$ denote the subset of $\Phi G$ consisting of
all $(A) \leq (K) $ in $\Phi G$. Lemmas \ref{openpoint} and
\ref{abelianfamily} imply that if $K$ is in $\aaa_n$, then
$\mathcal{S} (K)$ is a finite subset of $\aaa_n$ consisting of
open-closed points in $\Phi G$. We prove the proposition by
induction on the length of chains (totally ordered subsets) in the
partially ordered set $\mathcal{S} ( K ) $. If $(K)$ is minimal in
$\aaa_n$, then $\phi (G/K) = |W_G K| e_{K}$ and the claim is true.
Assume the claim is true for all $(A)$ such that all chains in
$\mathcal{S} (A)$ have length $(l-1)$ or less. If all chains of
 subgroups in $\mathcal{S} (K)$ have  length $l$ or less, then
we get
\[ \frac{\g_n }{ {|W_G K |}} \,  \phi (G /K ) =  \g_n e_{K} +
\textstyle\sum\limits_{(A) \lneqq (K)} \g_n m_A e_{A} \] where $m_A
= | (G/K )^A|/ |W_G K|$ are integers (as explained earlier in the beginning
of section 2). By the induction hypothesis we get that $\g_n e_{K} $
is in $\phi A(G)$, and $\g_n e_{K} = \textstyle\sum\limits
 c_A \phi ( G/A)$ where the sum is over ${(A) \leq (K)}$.
\end{proof}

 We conclude that  there is a stable
 map   $\alpha_n \col  S_{G} \rightarrow S_{G} $
 whose  degree function $d (\alpha_n ) \in C (G) $ is
   $ \sum_{(A) \in \aaa_n }    \g_n e^{G}_A
$. The degree of $ ( \g_n - \alpha_n )^A $ is $ 0 $ for all $(A) \in
\aaa_n$. Hence if $E_G$ satisfies the $n$-induction property, then
$(\g_n - \alpha_n )E_{G}^{\ast} =0 $. So
 $\g_n E_{G}^{\ast} $ is equal to $
\alpha_n E_{G}^{\ast} $.

The element $[G/H] $ in $A(G)$ corresponds, via the isomorphism
$\chi$, to
\[  S_{G}  \stackrel{\tau }{\longrightarrow}  G/H_+
 \stackrel{c }{\longrightarrow}S_{G} \] in $\pi_{0}^{G}  ( S_{G} )$.
Here $\tau$ is the transfer map and $c$ is the collapse map. Hence
we have that \[ [G/H] = \ind_{H}^G [\ast] . \] The isomorphism class
of a point $[\ast ] \in A(H) $ corresponds to the identity map in
$\pi_{0}^{H} ( S_{H}^0 )$.

The next result is the Artin induction theorem.
\begin{thm}[Artin induction theorem]  \label{image}
Assume that $E_{G}$ is a ring spectrum satisfying the $n$-induction
property. Then the integer $\g_n$ times the unit element in
$E_{G}^{*}$ is in the image of the induction map \[ \textstyle\bigoplus\nolimits_{(A) \in
\aaa_n} \text{ind}_{A}^{G} \col \textstyle\bigoplus\nolimits_{(A) \in \aaa_n} E^{0}_{A}
\rightarrow E^{0}_{G} \] where the sum is over representatives for
each conjugacy class $(A) \in \aaa_n$.
\end{thm}
\begin{proof}  We have that $\g_n 1 = \alpha_n 1 $  in $E^{\alpha }_G$.
 The theorem   follows from  the Frobenius reciprocity law
\ref{frobenius} and Lemma \ref{idempotent}.
More precisely, let $f \col S_{G} \rightarrow E_{G} $ represent an
element in $E^{0}_{G} $. Then $ [G/H ] \cdot f = \text{ind}_{H}^{G}
[c \wedge f ] $, where $c \wedge f \col G/H_+ \rightarrow E_G $
represents an element in $E^{0}_{H} \cong E^{0}_{G}( G/H_+ ) $.
\end{proof}

As a consequence of the Artin induction theorem we can reconstruct
$E_G (X)$, rationally, from all the $E_A (X)$ with $(A) \in \aaa_n$
and the restriction and conjugation maps. To do this we need the
double coset formula for compact Lie groups. The double coset
formula was first proved by M. Feshbach \cite{fes1}. We follow the
presentation given in \cite[IV.6]{lms}. To state the double coset
formula it is convenient to express the induction and restriction
maps between $E_H$ and $E_K$ for subgroups $H$ and $K$ of $G$ by
maps in the $G$-stable homotopy category. Let $H \leq K$ be
subgroups of $G$. There is a collapse map
\[c_{H}^K \col   G/H_+ \rarr G/K_+ \] and a transfer map
\[ \tau_{H}^K \col  G/K_+ \rarr G/H_+  , \]   which  induce restriction
and induction maps \cite[p.~204]{lms}. Let $g$ be an element in $G$.
Right multiplication by $g$ induces an equivalence of $G$-manifolds
\[ \beta_g \col  G/H_+ \rarr G/(g^{-1} H g)_{ +} . \] Consider
$G/H$ as a left $K$-space. The space is a compact differentiable
$K$-manifold so it has finitely many orbit types \cite[5.9.1]{td}.
The orbit type of an element $x$ is the $K$-isomorphism class of the
homogeneous space $K x$.
  The
stabilizer of the element $gH$ is $ K \cap g H g^{-1}$. The left
$K$-quotients of subspaces of $G/H$ consisting of all points of a
fixed orbit type are manifolds \cite[IV.3.3]{bre}.
 These manifolds are called
 the orbit type manifolds  of $K \backslash G /H$.
 We decompose the double coset space $K \backslash G
/H$ as a disjoint union of the connected components $M_i$ of all the
 orbit type manifolds.
 We are now ready to state the double coset formula.
\begin{thm}[Double coset formula] \label{double}  Let $G$ be a compact Lie group and
$H $ and $ K$ be closed subgroups of $G$. Then we have
\[ \tau_{H}^G \circ c_{K}^G \simeq \textstyle\sum\nolimits_{M_i}
z(M_i ) \, \beta_{g} \circ c_{K \cap g H g^{-1} }^{gHg^{-1}} \circ
\tau_{K \cap g H g^{-1} }^K \] where the sum is over orbit-type
manifold components $M_i$ and $g \in G$ is a representative of each
$M_i$. The integer $z (M_i)$ is the internal Euler characteristic.
It is the Euler characteristic of the closure of $M_i$ in $K
\backslash G /H$ subtracting  the Euler characteristic of its boundary.
\end{thm} Recall  that the transfer map
\[ \tau_{H}^K \col  G/K_+ \rarr G/H_+  \] is trivial if the Weyl group
$ W_K H$ is infinite. In particular, we have the following
\cite[II.17]{fes1} \cite[IV.6.7]{lms}.

\begin{lem} \label{torus} Assume   $H =K$ is a maximal torus $T$ in $G$. Then the sum in
  the double coset formula, Theorem \ref{double}, simplifies to
the   sum over elements   $g\in G $ representing each $gT$ in the
Weyl group $W_G T$ of $T$. \end{lem}

Let $E_G$ be a $G$-equivariant ring spectrum that satisfy the
$n$-induction property, and let $M_G$ be a module over $E_G$. There
is a restriction map
\[ M^{\alpha}_G (X) \rightarrow \text{Eq} [
\textstyle\prod\nolimits_{A} M^{\res \, \alpha}_A (X)
\rightrightarrows \textstyle\prod\nolimits_{K , L, g} M^{\res \,
\alpha}_{K\cap gLg^{-1}} (X) ] .\] The first product in the
equalizer
 is over representatives for conjugacy classes of $\aaa_n$,
  and the second
product is over pairs $K ,L$ of these subgroup representatives and
over representatives $g$ for each of the orbit-type manifold
components of $ K \backslash G / L$. The maps in the equalizer are
the two restriction (and conjugation) maps.

\begin{defn} Let $r$ be an integer.
 We say that a pair of maps \[ f \col  A \rightarrow B
\ \ \text{and} \ \ g \col B \rarr A \] between abelian groups is an
$r$-isomorphism pair if $ f \circ g = r $ and $ g \circ f = r$. A
map is an $r$-isomorphism if it is a map belonging to an
$r$-isomorphism pair.
\end{defn}

The Artin induction theorem implies the following Artin restriction
theorem.
\begin{thm} \label{standard}  Let $E_G$ be a $G$-ring spectrum satisfying the
 $n$-induction property.
 Let $M_G$ be an $E_G$-module spectrum. Then there exists a
 homomorphism
\[ \psi \col  \text{Eq} [ \textstyle\prod\nolimits_{A}
 M^{\res \, \alpha}_A (X) \rightrightarrows
\textstyle\prod\nolimits_{ K, L, g} M^{\res \, \alpha}_{K\cap
gLg^{-1}} (X) ] \rarr M^{\alpha}_G (X) \] such that the restriction
map and $\psi$ form  a $\g_n$-isomorphism pair.
 The first product in the equalizer
 is over representatives for conjugacy classes of $\aaa_n$,
  and the second
product is over pairs $K ,L$ of these subgroup representatives and
over representatives $g$ for each of the orbit-type manifold
components of $ K \backslash G / L$.
\end{thm}

An analogous result holds for homology.

\begin{proof} The following argument is standard. Our proof is close
to \cite[2.1]{mc}. We prove the result in the following generality.
Consider an element $r \in E^{0}_G $ in the image of the induction
maps from $E^{0}_{H_i}$ for a set of subgroups $H_i$ of $G$.
 In our case $r$ = $\g_n$ and the subgroups are
representatives for the conjugacy classes $\aaa_n$ by the Artin
induction theorem \ref{image}.

Let $r = \textstyle\sum\nolimits_{i=1}^k \ind_{H_i}^G r_i$ where
$r_i \in E_{H_i}$. Define $\psi$ by setting
\[ \psi ( \textstyle\prod\nolimits_{H_i} m_{H_i} ) = \textstyle\sum\nolimits_{i=1}^k \ind_{H_i}^G ( r_{i} \,  m_{H_i}
) . \] We have that \[\psi \circ \text{res} ( m ) =
\textstyle\sum\nolimits_{i=1}^k \ind_{H_i}^G ( r_{i} \, \res_{H_i}^G
m ) = \textstyle\sum\nolimits_{i=1}^k \ind_{H_i}^G ( r_{i} ) m = r m
.
\] The second  equality follows from the Frobenius reciprocity law
\ref{frobenius}.

We now consider the projection of $ \text{res} \circ \psi ( \prod
m_{H_i} )$ to $M_{K}$. It is
\[ \textstyle\sum\nolimits_{i} \res_{K}^G \ind_{H_i}^G ( r_i m_{H_i}) .\]
By the double coset formula, Theorem \ref{double}, this equals
\[ \textstyle\sum\nolimits_{i} \textstyle\sum\nolimits_{  K g H_i} z_i \  \ind_{g H_i g^{-1} \cap  K }^{K}
\res_{g H_i g^{-1} \cap K }^{g H_i g^{-1}} \beta_g ( r_i \, m_{H_i}
)\] where for each $i$ the sum is over representatives $ K g H_i$ of
components of orbit-type manifolds of the double coset $ K
\backslash G / H_i$ and $z_i$ is an integer. By our assumptions we
have that
\[ \res_{g H_i g^{-1} \cap K}^{g H_i g^{-1} }\ m_{ g H_i g^{-1} }
= \res_{g H_i g^{-1} \cap K}^{K}\ m_{K} .
\]  So by Frobenius reciprocity we get
\[ \textstyle\sum\nolimits_{i} \left( \textstyle\sum\nolimits_{ K g H_i } z_i \
\ind_{g H_i g^{-1} \cap K }^{K} \res_{g H_i g^{-1} \cap K }^{g H_i
g^{-1}} \beta_g
 r_i \right ) m_{K} .\] This
equals
\[ \textstyle\sum\nolimits_{i}  ( \res_{K}^{G } \ind_{H_i}^G ( r_i )) m_K  =
 \res_{K}^G (r) m_K . \qedhere\]
\end{proof}
For a fixed compact Lie group $G$ both the restriction map and the map $\Psi
$ are natural in $M_G$ and $X$.

\section{Brauer induction} \label{Brauer}
We present an integral induction theorem for cohomology theories
satisfying the $n$-induction property. We first discuss some classes
of subgroups.

\begin{defn} \label{nhyper}
A subgroup $H$ of $G$ is $n$-hyper if it has finite Weyl group and
there is an extension
\[ 0 \rightarrow A \rightarrow H \rightarrow  P
\rightarrow 1 \] such that:
\begin{enumerate}
\item
 $P$ is a finite  $p$-group for some prime number $p$,
 \item
$A$ is an abelian subgroup of $G$, such that $\omega (A)$ is
topologically generated by $n$ or fewer elements, and $ | A
/A^{\circ} | $ is relatively prime to $p$.
\end{enumerate}
\end{defn}
\begin{lem} \label{lem:nsub}
Let $H$ be an $n$-hyper subgroup of $G$ (for the prime $p$) and let
$K$ be a subgroup of $H$. Then $K$ is an $n$-hyper subgroup of $G$
(for the prime $p$) if $K$ has finite Weyl group in $G$.
\end{lem} \begin{proof} Let $H$ be as in Definition \ref{nhyper}.
Let $Q$ be the image of $K$ in $P$ under the homomorphism $ H \rarr P$.  Then $Q$ is a $p$-group.  Let $B$ be the kernel of the homomorphism $ K \rarr Q$. Then $B$ is a subgroup of $A$.  The abelian group $B$ is a direct sum of subgroups $C$ and $D$, such that $D$ is a  finite $p$-group  and  $| C / C^{\circ} |$ is relatively prime to $p$ (see Lemma \ref{absplit}).
Lemma \ref{abelianfamily} implies that $\omega (C)$ is in $\aaa_n$ since $C$ is a subgroup of $A$, which is in $\aaa_n$. The subgroup $C$ is normal in $K$ and the quotient group  $K / C$ is an extension of $D $ by $Q$, hence a $p$-group.
\end{proof}

We next describe the idempotent elements in the Burnside ring $A(G)$
localized at a rational prime. First we need some definitions. A
group $H$ is said to be $p$-perfect if it does not have a nontrivial
(finite) quotient $p$-group. The maximal $p$-perfect subgroup $H'_p$
of $H$ is the preimage in $H$ of the maximal $p$-perfect subgroup of
the group of components $H/H^{\circ}$. Let $H$ be a subgroup of a
fixed compact Lie group $G$. Let $H_p$ denote the
 conjugacy class
 $\omega (H'_p)$
 of subgroups of $G$ with finite Weyl group associated with $H'_p$.
  Let $\Phi_p G$ denote  the
subspace of $\Phi G$ consisting of conjugacy classes of all
$p$-perfect subgroups of $G$ with finite Weyl group in $G$. Let
$N_{(p)}$ denote the largest factor of an integer $N$ that is
relatively prime to $p$. The next result is proved for finite groups
in \cite[7.8]{td}, and for compact Lie groups in \cite[3.4]{fol}.
\begin{thm} \label{localizedidempotent}
Let $G$ be a compact Lie group. Let $H$ be a $p$-perfect subgroup
such that $(H) \in \Phi_{p} G$ is an open-closed point in $\Phi_{p}
G $. Let $N$ be an integer such that $N  e_H$ is in $\phi A(G)$.
Then there exists an idempotent element $I_{H,p} \in C ( \Phi G ,
\mathbb{Z}) $ such that
  $N_{(p)} I_{H,p}  \in \phi A(G)  $, and $I_{H,p} $ evaluated at $(K)$ is
$ 1 $ if $ K_p = ( H)$ and zero otherwise. In particular, $I_{H,p} $
is an idempotent element in the localized ring $\phi A(G)_{(p)}$.
\end{thm}

  An abelian group $A$ is $p$-perfect if and only if
  $| A / A^{\circ} |$ is relatively
prime to $p$ by Lemma \ref{absplit}. We can apply Theorem
\ref{localizedidempotent} to $p$-perfect abelian groups in $\aaa_n$
with $N = \g_n$ by Lemma \ref{openpoint} and Proposition
\ref{idempotent}. Let $ I_{(p,n)} $ be $ (|G|_n)_{(p)}
\textstyle\sum\nolimits_{(A) } I_{A,p} $ where the sum is over all
$(A) \in \aaa_n$ such that $|A/ A^{\circ} |$ is relatively prime to
$p$. The element $I_{(p,n)} \in C(G)$ is in $\phi A(G)$. The
function $I_{(p,n)}$ has the value $(|G|_n)_{(p)}$ at each conjugacy
class $(H)$ of the form \[ 0 \rarr S \rarr H \rarr P \rarr 1
\] where $P$ is a $p$-group and $S$ is abelian with $| S/
S^{\circ}|$ relatively prime to $p$ and $\omega (S) \in \aaa_n$, and
 $I_{(p,n)}$ has the value
 0 at all other elements of $\Phi G$.  In particular, $I_{(p,n)}$
 has the value
$(|G|_n)_{(p)}$ at each $A \in \aaa_n$ by Lemma \ref{absplit}.
 The greatest common divisor of $ \{
(|G|_n)_{(p)} \}$, where $p$ runs over primes $p$ dividing $|G|_n$,
is 1. Hence there is a set of integers $z_p$ such that $
\textstyle\sum\nolimits_{p | \, |G|} z_p |G|_{(p)} = 1 $. Let $I_n$
be the element $ \textstyle\sum\nolimits_{p }
 z_p I_{(p,n)} $ where the sum is over primes $p$ dividing $|G|$.
The element $I_n$ is not an idempotent element in $C(G)$. The
function $I_n \col \Phi G \rarr \zz $ has the value 1 for all $(A)
\in \aaa_n $. Since $I_n$ is in the image of $\phi \col A(G) \rarr
C(G)$, there is a map $\beta_n \col S_{G} \rarr S_{G}$ so that
the degree $d ( \beta_n)$ is $I_n$. The degree of $ (1 - \beta_n
)^A$ is zero for all $(A) \in \aaa_n$. Assume $E_G$ is a
$G$-equivariant ring spectrum satisfying the $n$-induction property.
 We get that
\[  E_{G}^{\ast} ( S_{G} ) = \beta_n  \cdot E_{G}^{\ast} ( S_{G} ) . \]

\begin{lem} \label{decomp} The element $I_n \in A(G)$ can be written as
\[ I_n = \phi ( \textstyle\sum\nolimits_{H_i } k_{H_i} [G/ H_i ] ) \]
where the subgroups $ H_i$ are $n$-hyper subgroups of $G$ for primes
dividing $|G|_n$, and $k_{H_i} $ are integers.
\end{lem}
\begin{proof}  We know that $I_n$ can be written in the above
form for some subgroups $H_i$ of $G$ with finite Weyl groups. Let $H_j$ be a maximal subgroup of $G$ in the
sum describing $I_n$. The value of $I_n$ at $(H_j) $ is $ k_{H_j}|
W_G H_j |$. This value is nonzero; hence the maximal subgroups $H_j$ are
  $n$-hyper subgroups of $G$. By Lemma \ref{lem:nsub}  all  $H_i$ are   $n$-hyper subgroups of $G$.
\end{proof}

The next result is the Brauer induction theorem.
\begin{thm}[Brauer induction theorem] \label{prebrauer}  Assume $E_{G}$ is a ring spectrum
 satisfying  the $n$-induction
property. Then the unit element 1 in $E_{G}^{\ast}$ is in the image
of the induction map \[ \textstyle\bigoplus\nolimits_{ (H)} \text{ind}_{H}^{G} (X) \col
\textstyle\bigoplus\nolimits_{(H)} E^{0}_{H} \rightarrow E^{0}_{G}
\]  where the sum is over   $n$-hyper subgroups of $G$ for
primes $p$ dividing $|G|_n$.
\end{thm}
\begin{proof}
We get as in Theorem \ref{image} that the unit element 1 is in the
image of the induction maps from all the subgroups $H_i$ in Lemma
\ref{decomp}.
\end{proof}
As a consequence we get the following Brauer restriction theorem.

\begin{thm}[Brauer restriction theorem] \label{Brauerinduction}
Let $E_G$ be a $G$-equivariant ring spectrum satisfying the
$n$-induction property, and let $M_G$ be a module over $E_G$. Then
the restriction map
\[ M^{\alpha}_G (X) \rightarrow \text{Eq}
 [ \textstyle\prod\nolimits_{H} M^{\res \, \alpha}_H (X)
\rightrightarrows \textstyle\prod\nolimits_{K , L , g} M^{\res \,
\alpha}_{K \cap gLg^{-1}} (X)]
\] is an isomorphism. The first product in the equalizer
 is over representatives for conjugacy classes of  $n$-hyper
subgroups of $G$ for primes dividing $|G|_n$, and the second product
is over pairs $K ,L$ of these subgroup representatives and over
representatives $g$ for each of the orbit-type manifold components
of $ K \backslash G / L$.
\end{thm}
\begin{proof} This follows from the proof of
  Theorem \ref{standard} and the
  Brauer induction theorem \ref{prebrauer}.
\end{proof}

An analogous result holds for homology.

\section{Induction theory for Borel  cohomology}   \label{Borel}
Let $k$ be a nonequivariant spectrum.
  The Borel cohomology and Borel homology on the category of based
$G$-spaces are
\[ k^{\ast} (X\wedge_G EG_+  ) \ \text{and} \
k_{\ast} ( \Sigma^{\text{Ad} (G)} X \wedge_GEG_+ ) . \]
 The adjoint representation $\text{Ad} (G) $ of $G$ is the tangent
vectors pace at the unit element of $G$ with $G$-action induced by
the conjugation action by $G$ on itself. If $G$ is a finite group,
then $\text{Ad} (G) =0$. Borel homology and cohomology can be
extended to an $RO(G)$-graded cohomology theory defined on the
stable equivariant homotopy category. We follow Greenlees and May
\cite{gm}.

Let $M_G$ be a $G$-spectrum. The geometric completion of $M_G$ is
defined to be
\[ c( M_G ) = F ( EG_+ , M_G ), \]
where $F$ denotes the internal hom functor in the $G$-equivariant
stable homotopy category. Let $f (M_G ) $ be $ M_G \wedge EG_+$. The
Tate spectrum of $M_G$ is defined to be
\[ t( M_G ) = F ( EG_+ , M_G ) \wedge \tilde{E}G .\]
The space $\tilde{E}G $ is the cofiber of the collapse map $ EG_+
\rarr S^0 $.

If $E_G$ is a ring spectrum, then $ c(E_G)$ and $t(E_G )$ are also
ring spectra \cite[3.5]{gm}. More precisely, $ c(E_G)$ is an algebra
over $E_G$, and $t (E_G )$ is an algebra over $ c(E_G)$. The product
on the spectrum $f(E_G) $ is not unital in general; however, it is a
$c (E_G)$ module spectrum. Hence we have the following.
\begin{pro}  Let $E_G$ be a $G$-ring spectrum satisfying the
$n$-induction property. Then $ c(E_G)$ and $t(E_G )$ are
$G$-equivariant ring spectra satisfying the $n$-induction property.
The spectrum $f (E_G )$ is a $ c(E_G)$-module spectrum.
\end{pro}

 Let $i \col  \mathcal{U}^G \rightarrow \mathcal{U}
$ be the inclusion of the universe $\mathcal{U}^G $ into a complete
$G$-universe $ \mathcal{U}
 $. Let $i_{\ast} k$ denote the
$G$-spectrum obtained by building in suspensions by
$G$-representations. If $k$ is a ring spectrum, then $i_{\ast} k$ is
a $G$-equivariant ring spectrum. The following result is proved in
\cite[2.1,3.7]{gm}.
\begin{pro} Let $k$ be a spectrum.   Let $X$ be a naive
$G$-spectrum (indexed on $\uu^G$). Then we have isomorphisms
\[ (c(i_{\ast} k ))_{G}^n  ( X )
 \cong k^n ( X \wedge_G EG_+ ) \ \text{and} \ (f (i_{\ast} k ))^{G}_n (X)  \cong
 k_n ( \Sigma^{\text{Ad} (G) } X \wedge_G EG_+ ) \]
where $\text{Ad} (G)$ is the adjoint representation of $G$.
\end{pro}

We next show that $ c (k)$ has the induction property (see
Definition \ref{inductionproperty}) when $k$ is a complex oriented
spectrum. For every compact Lie group $G$, there is a finite
dimensional unitary faithful $G$-representation $V$. Hence $G$ is a
subgroup of the unitary group $U(V)$. Let the flag manifold $F $
associated with  the $G$-representation $V$ be
 $ U(V)/T $ where $T$ is some  fixed  maximal torus of $U(V)$.  Let
$G$ act on $F$ via the embedding $G \leq U(V)$. The following result
is well known \cite[2.6]{hkr}.
\begin{pro} \label{injective} Let  $k$  be a complex oriented  spectrum.
Then the map
\[ k^{\ast} ( BG_+ ) \rightarrow k^{\ast} ( F_+ \wedge_G EG_+  ), \]
induced by the collapse map $F \rarr *$, is injective.
\end{pro}
\begin{thm} \label{abid} Given a compact Lie group $G$ let $N$
be an integer so that $ \aaa_N = \aaa$. Let $k$ be a complex
oriented spectrum. Then there exists an integer $d$ so that the
Borel cohomology $G$-spectrum $ c (i_{\ast} k)$ satisfies the
following: $ (\g_N - \alpha_N )^d c (i_{\ast} k)_* = 0 $ and $ (1-
\beta_N)^d c (i_{\ast} k)_* =0$.
\end{thm}
The class $\alpha_N$ is defined after Proposition \ref{idempotent}
and the class $\beta_N$ is defined before Lemma \ref{decomp}.
\begin{proof} An easy induction gives the following:
Let $Y$ be a $d$-dimensional $G$-CW complex, and let $E_G$ be a
cohomology theory. Assume that an element $r \in A(G)$ kills the $
E_{G}^* $-cohomology of all the cells in $Y$. Then $r^d$ kills $
E_{G}^* (Y)$.

  By Proposition \ref{injective} it
suffices to show that $ r^d (i_{\ast} k)^* (U/T) =0 $ when $ r$ is $
\g_N - \alpha_N$ or $1- \beta_N$. We have that $U/T$ is a finite
$G$-CW complex with orbit types $G/ (G \cap g T g^{-1})$ for $g \in
G$. Since both $\g_N - \alpha_N$ and $1-\beta_N$ restricted to any
abelian group are  0 by Corollary \ref{abeliandegree}, we get that
they kill the $c (i_{\ast} k)$-cohomology of all the cells in $Y=
U/T$. We can take $d $ to be the dimension of the $G$-CW complex
$U/T$.
\end{proof}

The Propositions \ref{injective} and \ref{abid} give Artin and
Brauer restriction theorems for $c(k)$ and $t(k)$ where $k$ is a
complex oriented spectrum. We state the theorems only for $c(k)$
applied to naive $G$-spectra.
\begin{thm}  Let $k$ be a complex oriented cohomology theory.
 Then for any naive $G$-spectrum $X$  the
restriction map from \(k^{\ast } ( X \wedge_G EG_+ ) \) to the
equalizer of \[ \textstyle\prod\nolimits_{(A) \in \aaa} k^{\ast } (X \wedge_A
 EG_+ ) \rightrightarrows \textstyle\prod\nolimits_{(K) , (L) \in \aaa \text{ and } KgL \in K\backslash  G / L }
  k^{\ast} ( X \wedge_{K\cap gLg^{-1}} EG_+ )
\] is a natural isomorphism after inverting $|G|$.
The first product
 is over representatives for conjugacy classes of $\aaa$, and the second
product is over pairs $K ,L$ of these subgroup representatives and
over representatives $g$ for each of the orbit-type manifold
components of $ K \backslash G / L$.
\end{thm}
\begin{thm}  Let $k$ be a complex oriented cohomology theory.
 Then for any naive $G$-spectrum $X$  the
restriction map from \(k^{\ast } ( X \wedge_G EG_+ ) \) to the
equalizer of \[ \textstyle\prod\nolimits_{H } k^{\ast } (X \wedge_H
 EG_+ ) \rightrightarrows
 \textstyle\prod\nolimits_{K,L, g} k^{\ast} ( X \wedge_{K \cap
gLg^{-1}} EG_+ )
\] is a natural isomorphism. The first product
 is over representatives for conjugacy classes of hyper
subgroups of $G$ for primes dividing $|G|$, and the second product
is over pairs $K ,L$ of these subgroup representatives and over
representatives $g$ for each of the orbit-type manifold components
of $ K \backslash G / L$.
\end{thm}

\begin{rem}
See also \cite{fes2} and \cite[IV.6.10]{lms}.
\end{rem}
Since $f( i_{\ast} k)$ is a $c( i_{\ast} k)$ module spectrum we get
 Artin and Brauer restriction theorems for the
  Borel homology $k_{\ast } ( EG \wedge_G \Sigma^{\text{Ad} (G)}
X )$ when $k$ is complex oriented.

\begin{thm}  Let $k$ be a complex oriented cohomology theory.
Then for any naive $G$-spectrum $X$ the restriction map from $
k_{\ast } (\Sigma^{\text{Ad} (G)} X \wedge_G EG_+ )$ to the
equalizer of
\[  \textstyle\prod\nolimits_{A} k_{\ast }
(\Sigma^{\text{Ad} (A) } X \wedge_A EG_+ ) \rightrightarrows
\textstyle\prod\nolimits_{K , L , g} k_{\ast} ( \Sigma^{\text{Ad}
(K\cap gLg^{-1}) } X \wedge_{K\cap gLg^{-1}} EG_+ )
\] is a natural isomorphism after inverting $|G|$.
The first product in the equalizer
 is over representatives for conjugacy classes of $\aaa$, and the second
product is over pairs $K ,L$ of these subgroup representatives and
over representatives $g$ for each of the orbit-type manifold
components of $ K \backslash G / L$.
\end{thm}
\begin{thm}  Let $k$ be a complex oriented cohomology theory.
 Then for any naive $G$-spectrum $X$  the
restriction map from $ k_{\ast } (\Sigma^{\text{Ad} (G)} X \wedge_G
EG_+ )$ to the equalizer of
\[  \textstyle\prod\nolimits_{H}
k_{\ast } (\Sigma^{\text{Ad} (H) } X \wedge_H EG_+ )
\rightrightarrows \textstyle\prod\nolimits_{K , L , g} k_{\ast} (
\Sigma^{\text{Ad} (K\cap gLg^{-1}) } X \wedge_{K\cap gLg^{-1}} EG_+
)
\] is a natural isomorphism. The first
product
 is over representatives for conjugacy classes of hyper
subgroups of $G$ for primes dividing $|G|$, and the second product
is over pairs $K ,L$ of these subgroup representatives and over
representatives $g$ for each of the orbit-type manifold components
of $ K \backslash G / L$.
\end{thm}

It is immediate from the definition that the induction and
restriction maps in Borel cohomology are given by the transfer and
the collapse maps
\[ \tau \col  S_{G} \rarr G/H_+  \ \  \text{and} \ \ c \col  G/H_+ \rarr S_{G}
 \] as follows: \[ k^{\ast} ( X \wedge_H EG_+ ) \cong
 k^{\ast} ((G/H_+ \wedge X ) \wedge_G  EG_+ )
 \stackrel{\tau^{\ast}}{\rarr}
k^{\ast} ( X \wedge_G EG_+ )\] and
\[ k^{\ast} ( X  \wedge_G EG_+ ) \stackrel{c_{\ast}}{\rarr}
k^{\ast} ( (G/H_+ \wedge X ) \wedge_G EG_+ ) \cong k^{\ast} (
 X  \wedge_H EG_+ ) . \]

When $G$ is a finite group the induction map in Borel homology is
induced from the collapse map $c$ and the restriction map is induced
from the transfer map $\tau$. This is more complicated for compact
Lie groups. The Spanier-Whitehead dual $ S^{- V} \thom{\nu G/H}
$ of $G/H_+$ (in the statement of Lemma \ref{homology}) is
equivalent to $ G \ltimes_H S^{- L(H) } $, where $ \ltimes $ is the
half smash product \cite[XVI.4]{ala} and $L(H)$ is the
$H$-representation on the tangent space at $eH$ in $G/H$ induced by
the action $ h, g H \mapsto hg H$.
  By considering
\[ H \rarr G \rarr G/H,
\] we see that $ Ad (G) $ is isomorphic to $Ad  (H)
\oplus L(H) $ as $H$-representations. By properties of the half
smash product we have that \[\left( G \ltimes_H S^{ - L(H) } \right)
\wedge S^{Ad( G)} \cong G \ltimes_H S^{Ad( G) - L(H) } \cong G
\ltimes_H S^{Ad ( H) } \] and
\[ ( G  \ltimes_H S^{Ad ( H) }   \wedge X )  \wedge_G EG_+ \cong
( S^{Ad ( H) } \wedge X ) \wedge_H EG_+ . \] Combined with Lemma
\ref{homology} this give a description of induction and restriction
maps.

Let $G$ be a finite group. We state some results from \cite{hkr} to
show that under some hypothesis a local complex oriented cohomology
theory satisfies $n$-induction. Let $k^{\ast}$ be a local and
complete graded ring with residue field of characteristic $p>0$.
Assume that
\[ k^{0} (BG) \rarr  p^{-1}k^{0} (BG) \]
is injective and that the formal group law
of $k^{\ast }$ modulo the maximal ideal has height $n$.

In \cite{hkr} the authors show that the restriction maps from
\[ p^{-1} k^{\ast} (BG) \]
 into the product of all  $p^{-1} k^{\ast} (BA)$
for all $p$-groups $(A) \in \aaa_n $ are injective.

 Hence if
$\beta $ is an element in the Burnside ring $A(G)$ of $G$ such that
$\text{deg}( \beta^A)=0$ for all $A$ in $\aaa_n$, then we have that
\[ \beta k^{\ast} ( BG ) = 0 .\] So the Borel cohomology $k^{\ast} (
X \wedge_G EG_+ ) $ satisfies the $n$-induction property. The
authors also show that $k^{\ast} ( X \wedge_G EG_+ ) $
 does not satisfy the $(n-1)$-induction property.

\section{0-induction, singular Borel cohomology} \label{singular}
We consider induction theorems for ordinary singular Borel
cohomology.
 Let $M$ be a $G/G^{\circ}$-module.  The Borel cohomology of
 an unbased
$G$-space $X$ with coefficients in $M$ is the singular cohomology of the Borel construction of $X$, 
\[ H^{\ast} ( X \times_G EG ;M ) ,\] with local coefficients via
\[ \Pi ( X \times_G EG ) \rightarrow \Pi ( \ast \times_G EG )
 \rightarrow
\pi_1 ( \ast \times_G EG ) \cong G / G^{\circ} \] where $\Pi$ is the
fundamental groupoid. Since $BG$ is path connected the fundamental
groupoid of $BG$ is non-canonically isomorphic to the one-object
category $\pi_1 (BG) $, which is the group of components of $G$.
Borel cohomology is an equivariant cohomology theory.
\begin{lem} The Borel cohomology of an unbased $G$-space $X$
with coefficients in a $G/ G^{\circ}$-module $M$ is represented in
the stable $G$-equivariant homotopy category by the geometric
completion of an Eilenberg-Mac\,Lane spectrum $ H \tilde{M}$, where
$\tilde{M}$ is a Mackey functor so that $ \tilde{M} (G/1) $ is
isomorphic to $M$ as a $G$-module.
\end{lem} Pre-composing with the $G$-equivariant suspension spectrum functor   is implicit in the
statement that the Borel cohomology is represented.
\begin{proof}
 The
Borel cohomology is isomorphic to the cohomology of the cochain
complex of $G$-homomorphisms from the $G$-cellular complex of $X
\times EG$
 to the $G$-module $M$  \cite[3.H]{hat}.
Assume that $\tilde{M}$ is a Mackey functor so that $ \tilde{M}
(G/1) $ is isomorphic to $M$ as a $G$-module.
 By the cell complex description
 of Bredon cohomology we get that the Borel  cohomology group
 is isomorphic to the
$G$-Bredon cohomology of $X \times EG$ with coefficients in
$\tilde{M}$. It follows from \cite[6.1]{gm} that there exist  Mackey
functors of the requested form.
\end{proof}
The relation between the geometric completion of Eilenberg-Mac\,Lane
spectra and classical Borel cohomology theory is also treated in
\cite[§6,§7]{gm}. From now on we consider the Borel cohomology
defined on $G$-spectra.

 The restriction map on the zeroth coefficient groups of the
 geometric completion of $ HM$
 is described by   the following commutative diagram:
\[ \xymatrix{ H^0 (BG_+ ; M )  \ar[r]^(.63){\cong} \ar[d]_{\text{res}}
& M^G \ar[d] \\ H^{0} ( B \ast_+ ; M ) \ar[r]^(.63){\cong} & M . }
\]
We have that $\text{res} ( \beta ) \text{res} (m) = \text{deg}
(\beta ) \text{res} (m) $ for any $\beta \in A(G) $ and $m \in M^G$.
Since the restriction map is injective we get that $ \beta m=
\text{deg} (\beta ) m $. So if $\text{deg} (\beta ) = 0$, then $
\beta H^0 (BG_+ ; M ) =0$. We have that $ \text{deg} ( \beta ) =
\text{deg} ( \beta^{T} )$ for any torus $T$ in $G$ by Proposition
\ref{omega}. So singular Borel cohomology satisfies the 0-induction
property.
The argument above applies more generally to show that Bredon
homology and cohomology with Mackey functor coefficients $M$ such
that $ M (G/G ) \rarr M (G/e)$ is injective (or, alternatively,
$M(G/e) \rarr M(G/G)$ is surjective) satisfies 0-induction.

 The  Artin restriction  theorem  gives a refinement of
Borel's description of rational Borel cohomology. Recall Lemma \ref{torus}. See also
\cite[II.3]{fes2}.
\begin{thm}
Let $G$ be a compact Lie group, $X$ a $G$-spectrum, and $M$ a
$G/G^{\circ}$-module. Then the restriction map
\[ H^{\ast} ( X \wedge_G EG_+ ;M ) \rightarrow H^{\ast} ( X
\wedge_{T} EG_+ ; M )^{W_G T} \] is a $ |W_G T |$-isomorphism.
\end{thm}
In particular, with $X= S^0$ we get that
\[ H^{\ast} ( BG_+ ;M ) \rightarrow H^{\ast} ( B T_+ ; M )^{W_G T} \]
is a $ |W_G T |$-isomorphism.  When $G$ is a finite group this says that  $ H^0 ( BG_+ ; M) \rarr M^G$ is a $|G|$-isomorphism and $|G| \, H^{k} ( BG_+ ;M ) = 0 $, for $k > 0$.

 We next give the  Brauer restriction  theorem.

\begin{thm} Let $G$ be a
compact Lie group, $X$ a $G$-spectrum and $M$ a $G/ G^{\circ}
$-module.
 Fix a maximal torus $T$ in $G$. Then the restriction map
\[ H^{\ast} ( X \wedge_G   EG_+ ;M ) \rightarrow
\lim_K H^{\ast} ( X \wedge_K EG_+ ; M ) \] is an isomorphism. The
limit is over all  subgroups $K$ of $G$ with finite Weyl group that have a normal
abelian
 subgroup $A$ of $K$ such that  $A \leq T$ and   $K/A$ is a
$p$-group for some prime $p$ dividing $|W_G T|$. The maps in the
limit are restriction maps and conjugation maps.
\end{thm}

  We use singular homology  with local coefficients.
\begin{thm}
Let $G$ be a compact Lie group, $X$ a $G$-spectrum, and $M$ a $G/
G^{\circ}$-module. Then the restriction map (induced by t)
\[ H_{\ast} (  S^{\text{Ad} (G)} X \wedge_G EG_+ ;M )
\rightarrow H_{\ast} ( S^{\text{Ad} (T)} X \wedge_{T} EG_+ ; M
)^{W_G T}
\] is a $ |W_G T |$-isomorphism.
\end{thm}   There is also a
Brauer restriction theorem for homology.
\begin{thm} Let $G$ be a
compact Lie group, $X$ a $G$-spectrum and $M$ a $G/ G^{\circ}
$-module.
 Fix a maximal torus $T$ in $G$. Then the restriction map
\[ H_{\ast} ( S^{\text{Ad} (G)} X \wedge_G   EG_+ ;M ) \rightarrow
\lim_K H_{\ast} ( S^{\text{Ad} (K)} X \wedge_K EG_+ ; M ) \] is an
isomorphism. The limit is over all the subgroups $K$ of $G$ with finite Weyl group that
have a normal abelian
 subgroup $A$ of $K$ such
  that $A \leq T$ and   $K/A$ is a
$p$-group for some prime $p$ dividing $|W_G T|$. The maps in the
limit are restriction   and conjugation maps.
\end{thm}

\section{1-induction, $K$-theory} \label{ktheory}
In this section we consider equivariant $K$-theory $K_G$. For
details on $K_G$, see \cite{seg2}. An element $g \in G$ is said to be
regular if the closure of the cyclic subgroup generated by $g$ has
finite Weyl group in $G$.
\begin{defn}
Let $\rho G$ denote the space of conjugacy classes of regular
elements in $G$. Define $r \col \rho G \rightarrow \Phi G $ by
sending a regular element $g$ in $G$ to the closure of the cyclic
subgroup generated by $g$.
\end{defn}
The space $\rho G $ is given the quotient topology from the subspace
of regular elements in $G$. The map $r$ is continuous since two
regular elements in the same component of $G$
 generate conjugate cyclic  subgroups  \cite[1.3]{seg1}.
Let $C (X , R)$ denote the ring of continuous functions from a space
$X$ into a topological ring $R$. The map $r$, together with the
inclusion of $\mathbb{Z} $ in $\mathbb{C}$,
 induces a  ring homomorphism
\[ r^{\ast} \col  C( \Phi G , \mathbb{Z} )
\rightarrow C( \rho G , \mathbb{C} ). \] The ring of class functions
on $G$ is a subring of $C( \rho G , \mathbb{C} )$ since the regular
elements of $G$ are dense in $G$. Let $R(G)$ denote the (complex)
representation ring of $G$. Let $\chi \col R(G) \rightarrow C( \rho
G , \mathbb{C} )$ be the character map.
 The map  $\chi $ is an
injective ring map. Let $V$ be a $G$-representation. The value of
$\chi (V)$ at (a regular element) $g \in G$ is the trace of $ g\col
V \rightarrow V $. This only depends on the isomorphism class of
$V$.

 We now give a description of the induction map $
\ind_{H}^G \col R(H) \rightarrow R(G) $ for the coefficient ring of
equivariant $K$-theory
   \cite{oli,seg1}.
  Let $\xi$ be an $H$-character.
On a regular element $g$ in $G$ the induction map is given by
\[ \ind_{H}^G \xi (g)  = \textstyle\sum\nolimits_{kH} \xi (k^{-1} g k ) , \]
where the sum is over the finite fixed set $(G / H)^{g }$.
If $x \in R(H)$, then $ \ind_{H}^G \xi (x) $ is in the image of $
\xi \col R(G) \rarr C ( \rho G , \mathbb{C} ) $ \cite[2.5]{oli}. So
we get a well defined map $\ind_{H}^G \col R(H) \rarr R(G)$. This agrees with
the induction map defined in section \ref{cohomology} by
\cite[5.2]{nis} and \cite[§2]{seg1}.

The unit map in $G$-equivariant $K$-theory induces a map $ \sigma
\col A(G) \rightarrow R(G) $. It is the generalized permutation
representation map
\[  \sigma ( [X] )  = \textstyle\sum\nolimits_{i} (-1)^i  [ H^i
(X; \mathbb{C}) ] , \] where $[ H^i (X; \mathbb{C} ) ] $ is the
isomorphism class of the $G$-representation
 $ H^i (X; \mathbb{C} ) $ \cite[7]{tdb}.
As pointed out in \cite[7]{tdb}, see also \cite[5.3.11]{td}, the
following diagram commutes:
\[ \xymatrix{ A(G) \ar[r]^(.4){\phi}  \ar[d]^{\sigma}   &
C( \Phi G , \mathbb{Z} ) \ar[d]^{r^{\ast}} \\
R(G) \ar[r]^(.4){\chi} & C( \rho G , \mathbb{C} ) . } \]

\begin{lem} \label{oneind}
Equivariant $K$-theory $K_G$ satisfies the $1$-induction
 property. \end{lem} \begin{proof}
Since $\chi$ is injective it suffices to note that the map $r \col
\rho G \rarr \Phi G$ factors through $\aaa_1$.
\end{proof}
The Artin restriction theorem for equivariant $K$-theory follows
from Theorem \ref{standard}.
\begin{thm} For every $G$-spectrum $X$ the restriction map
 \[ K_G
( X) \rarr \text{Eq}[ \textstyle\prod\nolimits_{A} K_A (X)
\rightrightarrows \textstyle\prod\nolimits_{H ,L , HgL} K_{ H \cap g L
g^{-1} } (X) ]\] is a $|G|_1$-isomorphism. The first product in the
equalizer is over representatives for conjugacy classes of
topologically cyclic subgroups of $G$, and the second product is
over pairs of these subgroups and elements $HgL \in H\backslash G / L$. \end{thm} It suffices to pick a representative for each path
component of the $H$-orbit space of each submanifold of $ G /L$
consisting of points with a fixed orbit type under the $H$-action.

 When $G$ is connected, the maximal torus $T$ is  the only
conjugacy class of subgroups of $G$ with a dense subgroup generated
by one element and with finite Weyl group by Example \ref{cartan}.
So when $G$ is connected the restriction map
\[ K_G (X) \rightarrow K_T (X)^{W_G T }
\]  is a
$|W_G T|$-isomorphism.

We give an explicit description of $R(G)$ up to $\g_1$-isomorphism
using the Artin restriction theorem. Let $A^{\ast}=\hom ( A , S^1 )$
denote the Pontrjagin dual of $A$. The elements of $A^{\ast}$ are
the one dimensional unitary representations of $A$. All irreducible
complex representations of a compact abelian Lie group are one
dimensional. We verify that when $A$ is a compact abelian Lie group
the canonical map
\[ \zz [ A^{\ast} ] \rarr R (A) \] is an isomorphism.
A subgroup inclusion $f \col H \rarr L$ induces a restriction map $
f^{\ast} \col L^{\ast} \rarr H^{\ast} $ of representations. By the
Artin restriction theorem  there is an (injective) $|W_G
T|$-isomorphism
\[ R(G)
 \rarr \text{Eq}[
\textstyle\prod\nolimits_{A} \zz [ A^{\ast} ] \rightrightarrows
\textstyle\prod\nolimits_{H ,L , HgL} \zz [({ H \cap g L g^{-1} }
)^{\ast} ] \, ] ,\]
where $A, H$, and $L$ are  representatives for conjugacy classes of
topologically cyclic subgroups of $G$, and $HgL$ are representatives for each of the orbit-type manifold components
of    the coset
$H\backslash G / L$.

The Brauer induction theorem \ref{prebrauer} and Lemma \ref{oneind}
applied to $K_G$ give that the identity element in $R(G)$ can be
induced up from the 1-hyper subgroups of $G$. This was first proved
by G.~Segal \cite[3.11]{seg1}. We get the following Brauer
restriction theorem from Theorem \ref{Brauerinduction} and Lemma
\ref{oneind}.
\begin{thm}  For every $G$-spectrum $X$ the restriction map
\[ K_G ( X) \rarr \text{Eq}[ \textstyle\prod\nolimits_{H} K_H (X)
\rightrightarrows \textstyle\prod\nolimits_{H ,L,g } K_{ H \cap g L
g^{-1} } (X) ]\] is an isomorphism. The first product in the
equalizer is over representatives for conjugacy classes of 1-hyper
subgroups of $G$, and the second product is over pairs of these
subgroups and representatives $g$ for each of the orbit-type manifold components
of the coset $H\backslash G /L$.
\end{thm}

There is another equivariant orthogonal $K$-theory $KO_G$ obtained by
using real instead of complex $G$-bundles. It is a ring spectrum and
the unit map of $K_G$ factors as
\[ S_{G} \rarr KO_G \rarr K_G , \] where the last map is induced by
tensoring the real $G$-bundles by $\mathbb{C}$. The coefficient ring
of $KO_G$ is $RO(G)$ in degree 0. Since $RO (G)$ injects into $R(G)$
we get that $KO_G$ satisfies the 1-induction property. Another
example of a ring spectrum that satisfies 1-induction is Greenlees'
equivariant connective $K$-theory \cite{gre}.

The restriction map $R(G) \rarr \prod R(C)$ is injective for all
compact Lie groups $G$, where the product is over topological cyclic
subgroups of $G$. When $G$ is connected we
 even have that $ R(G) \rarr R(T)^{W_G T}$ is an isomorphism.
Sometimes this can be used to prove stronger induction and detection
type results about $K_G$ than the results we get from our theory
\cite{jac,mc}. See also Proposition 3.3 in \cite{ati}. An
F-isomorphism theorem for equivariant $K$-theory has been proved by
Bojanowska \cite{boj}.

\begin{rem} It is reasonable to ask if tom Dieck's equivariant
complex cobordism spectrum satisfies the induction property
\cite{tdc}. Brun has shown that this is not the case. In fact
the unit map $A(G) \rarr MU_{G}^0$ is injective when $G$ is a finite
group \cite{bru}.
\end{rem}

\end{document}